\documentclass[12pt,a4paper,twoside]{article}
 \usepackage[latin1]{inputenc}
\usepackage[T1]{fontenc}
\usepackage[margin=2.5cm]{geometry}
\usepackage{lmodern}
\usepackage[french,english]{babel}
\usepackage{bef_aras_nothm}
\usepackage{amsmath,amsthm,amssymb}
\usepackage{fixmath}
\usepackage{upgreek}
\usepackage{enumitem}
\usepackage{amsfonts}
\usepackage{microtype}
\usepackage{float}
\usepackage{mathrsfs}

\usepackage{hyperref}
\hypersetup{
     colorlinks   = true,
     citecolor    = blue
}
\makeatletter
\newcommand{\oset}[3][0ex]{%
    \mathrel{\mathop{#3}\limits^{
            \vbox to#1{\kern-2\ex@
                \hbox{$\scriptstyle#2$}\vss}}}}
\makeatother

\usepackage{verbatim}

\usepackage[headsepline]{scrlayer-scrpage}

\clearscrheadfoot
\ohead{\pagemark}
\ihead{\headmark}
\automark[subsection]{section}

\theoremstyle{plain}

\allowdisplaybreaks[1]

\numberwithin{equation}{section}

\newtheorem{thm}{Theorem}[section]
\newtheorem{lem}[thm]{Lemma}

\theoremstyle{definition}

\begin{document}

\newpage
\setcounter{page}{1}
\clearscrheadfoot
\ohead{\pagemark}
\ihead{\headmark}

\title{A generalization of the Moreau--Yosida regularization}

\author{Aras Bacho\footnotemark[2]}

\date{}
\maketitle

\footnotetext[1]{Ludwig-Maximilians-Universit\"{a}t M\"{u}nchen, Mathematisches Institut, Theresienstr. 39, 80333 M\"{u}nchen, Germany.}

\begin{abstract}
In many applications, one deals with nonsmooth functions, e.g., in nonsmooth dynamical systems, nonsmooth mechanics, or nonsmooth optimization. In order to establish theoretical results, it is often beneficial to regularize the nonsmooth functions in an intermediate step. In this work, we investigate the properties of a generalization of the \textsc{Moreau--Yosida} regularization on a normed space where we replace the quadratic kernel in the infimal convolution with a more general function. More precisely, for a function $f:X \rightarrow (-\infty,+\infty]$ defined on a normed space $(X,\Vert \cdot \Vert)$ and given parameters $p>1$ and $\varepsilon>0$, we investigate the properties of the generalized \textsc{Moreau--Yosida} regularization given by 
\begin{align*}
f_\varepsilon(u)=\inf_{v\in X}\left\lbrace \frac{1}{p\varepsilon} \Vert u-v\Vert^p+f(v)\right\rbrace \quad ,u\in X.
\end{align*} We show that the generalized \textsc{Moreau--Yosida} regularization satisfies the same properties as in the classical case for $p=2$, provided that $X$ is not a \textsc{Hilbert} space. We further establish a convergence result in the sense of \textsc {Mosco}-convergence as the regularization parameter $\varepsilon$ tends to zero. 
\end{abstract} \textbf{Keywords} \textsc{Moreau--Yosida} regularization $ \cdot $ Convex analysis $ \cdot $ p-duality map $ \cdot $ \textsc{G\^{a}teaux} differentiability $ \cdot $ \textsc{Mosco}-convergence $ \cdot $  Nonsmooth analysis
\section{Introduction}
\subsection{Preliminaries and notation} 
We denote $(X,\Vert \cdot \Vert)$ a normed space and $(X^*,\Vert \cdot \Vert_*)$ its topological dual space. The duality pairing between $X^*$ and $X$ is denoted by $\langle \cdot, \cdot \rangle$. For a functional $f:X\rightarrow (-\infty,+\infty]$, the effective domain is defined by $\DOM(f):=\lbrace u\in X: f(u)<+\infty\rbrace$. The function $f$ is called proper if $\DOM(f)\neq \emptyset$. Furthermore, the subdifferential of $f$ in the sense of convex analysis is given by
\begin{align*}
    \partial f(u):=\lbrace \xi \in X^* : f(u)-f(v)\leq \langle \xi, v\rangle\, \text{ for all }v\in X \rbrace.
\end{align*} The functional $f$ is called subdifferentiable in $u\in X$ if $\partial f(u)\neq \emptyset.$ The domain of the subdifferential is defined by $\DOM(\partial f):=\lbrace u\in X: \partial f(u)\neq \emptyset \rbrace$. For a general proper, lower semicontinuous, and convex functional $f:X \rightarrow (-\infty,+\infty]$ on a normed space $(X,\Vert \cdot \Vert)$, the classical \textit{Moreau--Yosida} \textit{regularization} of $f$ is defined via
\begin{align} \label{class.MYR}
f_\varepsilon(u)=\inf_{v\in X}\left\lbrace \frac{1}{2\varepsilon} \Vert u-v\Vert^2+f(v)\right\rbrace \quad ,u\in X,
\end{align} where $\varepsilon>0$ is called the \textit{regularization parameter}. It is well known that the geometrical properties of the dual space $X^*$ are intimately related to the regularity properties of the regularization $f_\varepsilon$, see, e.g., \textsc{Barbu} \cite{Barb10NDMT} and \textsc{Barbu \& Precupanu} \cite{BarPre86COBS}. Roughly speaking, the better the geometrical properties of the dual space $X^*$ are, the better the regularization becomes. In the following, we elaborate on this in more detail. To do so, we recall the definition of the duality map $F_X:X \rightrightarrows X^*$, which is given by the set $F_X(v):=\lbrace \xi \in X^* :\langle \xi,v\rangle=\Vert v\Vert^2=\Vert \xi \Vert_*^2\rbrace$. It is well known that the duality map is given by the subdifferential of the mapping $u\mapsto \frac{1}{2}\Vert u\Vert^2$, i.e., $F_X(u)=\partial(\frac{1}{2}\Vert u\Vert^2)$ for all $u\in X$. Furthermore, it is easily checked that for all $u\in X$, the set $F_X(u)$ is non-empty, convex, bounded, and weak$*$-closed\footnote{Therefore, $F_X(u)$ is weak$*$-compact.}, see, e.g., \textsc{Barbu \& Precupanu} \cite[Section 1.2.4]{BarPre86COBS}. The duality map also has a geometrical interpretation: by the \textsc{Hahn--Banach} theorem, see, e.g., \textsc{Br\'{e}zis} \cite[Theorem 1.1, p. 1]{Brez11FASS}, for $u\in X$, there holds 
\begin{align*}
\Vert u\Vert=\max_{\overset{\zeta \in X^*}{\Vert \zeta\Vert_*=1}}\langle \zeta,u\rangle=\max_{ \overset{\zeta \in X^*}{\Vert \zeta\Vert_*=\Vert u\Vert}}\frac{\langle \zeta,u\rangle}{\Vert u \Vert} \geq \frac{\langle \xi,u\rangle}{\Vert u \Vert} \quad \text{for all $\xi\in X$ with $\Vert \xi \Vert_*=\Vert u\Vert$}.
\end{align*} Thus, an element of the dual space belongs to the duality map $\xi^*\in F_X(u)$ if and only if it solves the maximization problem
\begin{align}\label{DM.geometric}
\max_{ \overset{\zeta \in X^*}{\Vert \zeta\Vert_*=\Vert u\Vert}}\frac{\langle \zeta,u\rangle}{\Vert u \Vert},
\end{align} for which the set of maximizers is non-empty. In other words, $\xi^*$ generates a closed supporting hyperplane to the closed ball $\overline{B}(0,\Vert u\Vert)$.\\\\ 
Furthermore, we call a norm \textit{smooth} if and only if the duality map is single-valued, or geometrically speaking, each supporting hyperplane which passes through a boundary point of the sphere $S(0,\Vert u\Vert)$ with radius $\Vert u\Vert$ is also a tangential hyperplane. We call a normed space smooth if there is an equivalent smooth norm. From \eqref{DM.geometric}, it is then readily seen that if the dual space $X^*$ is strictly convex, i.e., the dual norm $\Vert\cdot\Vert_*$ is strictly convex, the element which generates the supporting hyperplane is unique, meaning that the duality map $F_X(u)$ is single-valued. In this case, the duality map is also demicontinuous \footnote{A map $f: X \rightarrow Y$ between two normed spaces $X$ and $Y$ is called demicontinuous if it is strong-to-weak* continuous.}, which implies that the norm on $X$ is \textsc{G\^{a}teaux} differentiable. If the dual space $X^*$ is uniformly convex\footnote{The normed space $X$ is called uniformly convex if for every $0<\varepsilon\leq 2$ there exists $\delta>0$ such that for any two vectors $x,y\in X $ with $\Vert x\Vert=\Vert y\Vert=1$ the condition $\Vert x-y\Vert \geq \varepsilon$ implies that $\left \Vert \frac{x+y}{2}\right \Vert\leq 1-\delta$. An uniformly convex space is in particular strictly convex.}, then the duality map is uniformly continuous on every bounded subset of $X$ and the norm on $X$ is uniformly \textsc{Fr\'{e}chet} differentiable in the sense that the limit
\begin{align*}
\lim_{\lambda\rightarrow 0}\frac{\Vert u+\lambda v\Vert-1}{\lambda}
\end{align*} exists uniformly in $x,y\in S(0,1)$, see \cite{BuiT02TNDM,Barb10NDMT}. Obviously, the regularity of the norm of a \textsc{Banach} space is deeply related to the geometrical properties of its dual space. If $X$ is a reflexive \textsc{Banach} space, then by the renorming theorem due to \textsc{Asplund} \cite{Asp67AVNO}, there always exist equivalent norms of $X$ and the dual space $X^*$ such that both $X$ and $X^*$ equipped with these norms are strictly convex and smooth, see \textsc{Barbu \& Precupanu} \cite[Theorem 1.105, p. 36]{BarPre86COBS}. Consequently, a reflexive \textsc{Banach} space can be equipped with an equivalent \textsc{G\^{a}teaux} differentiable norm such that the duality map is demicontinuous. It is well-known that a \textsc{Hilbert} space, in particular, is reflexive and that the duality map is identical with the \textsc{Riesz} isomorphism between the \textsc{Hilbert} space and its dual. For a more detailed discussion about the geometry of \textsc{Banach} spaces, and in particular with regard to the duality maps, we refer the interested reader to \cite{Barb76NSDE,Barb10NDMT,BarPre86COBS,BuiT02TNDM,Bynu71COUC,Bynu76WPLB,Dies75GBST,Zeme91SMLC}.

\subsection{Literature review}
The classical \textsc{Moreau--Yosida} regularization as defined in \eqref{class.MYR} has been studied extensively and has been employed successfully in many applications in order to circumvent the lack of regularity. The properties for the classical \textsc{Moreau--Yosida} regularization can, for reflexive \textsc{Banach} spaces, be found in, e.g., \textsc{Barbu} \cite{Barb10NDMT} and \textsc{Barbu \& Precupanu} \cite{BarPre86COBS} and for \textsc{Hilbert} spaces in, e.g., Attouch \cite{Atto84VCFO} and \textsc{Moreau} \cite{More65PDEH}. More general infimal convolutions defined by 
\begin{align}\label{infimal}
    (f\square g)(u):=\inf_{v\in X}\lbrace f(u-v)+g(v)\rbrace
\end{align} for proper, lower semicontinuous and convex functionals $g$ and $f$ defined on a \textsc{Hilbert} space has been studied in \textsc{Bauschke \& Combettes} \cite{BauCom17CAMH}. In particular, the \textsc{Pasch--Hausdorff} envelope, i.e., $g(v)=\beta \Vert v\Vert$, the \textsc{Moreau} envelope, i.e., $g(v)=\frac{1}{\gamma 2} \Vert v\Vert^2$, and the case $g(v)=\frac{1}{\gamma p} \Vert v\Vert^p, p>1,$ have been studied. 
The present work generalizes the previous works by showing that the \textsc{Moerau--Yosida} regularization, henceforth called p-\textsc{Moerau--Yosida} regularization, for the kernel $g(v)=\frac{1}{\gamma p} \Vert v\Vert^p, p>1$ defined on a reflexive \textsc{Banach} space satisfies all the properties as the classical \textsc{Moerau--Yosida} regularization \textsc{G\^{a}teaux} differentiability except from the \textsc{Lipschitz} continuity of the \textsc{G\^{a}teaux} derivative of the regularization in the case the underlying space $X$ is a \textsc{Hilbert} space. In addition, we show the convergence of the p-\text{Moerau--Yosida} regularization in the sense of \textsc{Mosco} as the regularization parameter vanishes.\\\\
A crucial assumption in all the previous results is the convexity of the functional $f$. It is remarkable that similar results have been obtained for non-convex functionals $f$ that are defined on a \textsc{Hilbert} space via the the so-called \textsc{Lions--Lasry} regularization introduced by P.L. \textsc{Lions} and \textsc{Lasry} \cite{LasLio86RRHS}. The \textsc{Lions--Lasry} regularization of a proper function $f:H\rightarrow (-\infty,+\infty]$ that is minorized by a quadratic function is defined by 
\begin{align}\label{Lions.Lasry}
    (f_\lambda)^\mu(u):=\sup_{v\in H} \inf_{w\in H} \left \lbrace f(w)+ \frac{1}{2\lambda}\Vert v-w\Vert^2-\frac{1}{2\mu} \Vert u-v\Vert^2\right \rbrace.
\end{align}  Similarly, one can define a regularization for a proper function $g$ that is majorized by a quadratic function. Among other properties, it has been shown in \textsc{Attouch} and \textsc{Aze} \cite{AttAze93ARFH} that these functions are \textsc{Fr\'{e}chet} differentiable with \textsc{Lipschitz} continuous derivative, weakly- or $\lambda$-convex (i.e. convex up to a square), satisfy $(f_\lambda)^\mu\leq f$, and that $(f_\lambda)^\mu(u)$ coincides with the \textsc{Moreau--Yosida} regularization when $f$ is convex. The \textsc{Fr\'{e}chet} differentiability has also been shown for a more general class of kernels which includes the  class of \textsc{Young} functions. For quadratic kernels, these results have been partially extended by \textsc{Str\"{o}mberg} \cite{Stro96ORBS} to the case where $X$ is a Banach space whose norm and dual norm are (locally) uniformly rotund, i.e., if $\Vert \cdot \Vert^2$ and $\Vert \cdot \Vert_*^2$ are (locally) uniformly convex functions. \textsc{Penot} \cite{Peno98PM} has studied the \textsc{Fr\'{e}chet} differentiability of the infimal convolution \eqref{infimal} in relation to the proximal mapping  
\begin{align*}
    P_{f,g}(u):=\lbrace v\in X: (f\square g)(u)= f(u-v)+g(v)\rbrace.
\end{align*} In particular, it has been shown that the non-emptiness of $P_{f,g}(u)$ is related to certain properties of the (\textsc{Fr\'{e}chet} or \textsc{Hadamard}) subdifferential of $f$, or under smoothness assumptions on  the \textsc{Banach} space $X$, the \textsc{Fr\'{e}chet} derivative of $f$.  \textsc{Penot} and \textsc{Ngai} \cite{VanPen16SORF} have further extended the result by imposing a milder growth condition on the function $f$. In addition, the authors studied the \textsc{Lions--Lasry} regularization for more general kernels, i.e., 
\begin{align*}
    (f_\lambda)^\mu(u):=\sup_{v\in X} \inf_{w\in X}\left \lbrace f(w)+ \frac{1}{\lambda}g(\Vert v-w\Vert)-\frac{1}{\mu} g(\Vert u-v\Vert)\right \rbrace
\end{align*} for a convex, monotonically increasing, coercive, and  continuously differentiable function $h:\mathbb{R}^+\rightarrow \mathbb{R}^+$ with $h(0)=0$. \\\\ However, for non-convex functionals $f$, none of the 
previous results duplicate our results, and it is subject to future work to reproduce our results for the non-convex case with the aid of the previous  results for non-convex functionals.  We refer the interested reader to \cite{Bern10LLRI,BeThZl11PEUS, JoThZa14DPME} and the references therein for more results in the non-convex case. The references presented here are indeed not exhaustive.

\section{Main result}
The question arises: if and to what extent the properties of the duality map are related to the regularization properties of the \textsc{Moreau--Yosida} regularization. We will see that the properties of the duality map are inherited by the subdifferential of the \textsc{Moreau--Yosida} regularization. In fact, we will answer the question for the more general so-called $p$-\textsc{Moreau--Yosida} regularization, which for $p>1$, is given by  
 \begin{align}\label{def:PMYR}
f_\varepsilon(u)=\inf_{v\in X}\left\lbrace \frac{\varepsilon}{p} \left \Vert \frac{u-v}{\varepsilon}\right \Vert^p+f(v)\right\rbrace \quad ,u\in X.
\end{align} 
The reason why we want to study $p$-\textsc{Moreau--Yosida} regularization is simply because it maintains the growth of the functional $f$ if it has $p$-growth, see, e.g., \cite{Bach20DNEI, Bach21ONSA}.\\\\
The following lemma shows some basic properties of the $p$-\textsc{Moreau--Yosida} regularization on general normed spaces. 
\begin{lem} \label{le:MYR.basic}  Let $f:X\rightarrow (-\infty,+\infty]$ be a proper and convex functional, and, for $\varepsilon>0$ and $p>1$, let $f_\varepsilon$ be the $p$-\textsc{Moreau--Yosida} regularization defined by \eqref{def:PMYR}. Then, $f_\varepsilon: X \rightarrow \mathbb{R}$ is finite, convex, and locally \textsc{Lipschitz} continuous. If, in addition, $f$ is  lower semicontinuous and $X$ is a reflexive \textsc{Banach} space, then the infimum in $f_\varepsilon(u)=\inf_{v\in X}\left\lbrace \frac{\varepsilon}{p} \left \Vert \frac{u-v}{\varepsilon}\right \Vert^p+f(v)\right\rbrace$ is attained at every point $u\in X$.
\end{lem}
\begin{proof}[Proof] Let $\tilde{u} \in \DOM(f)\neq \emptyset$. Then, on the one hand, there holds
\begin{align}\label{eq:bounded}
f_\varepsilon(u)\leq \frac{1}{p\varepsilon^{p-1}}\Vert u-\tilde{u}\Vert^p +f(\tilde{u})<\infty \quad \text{for every }u\in X.
\end{align} On the other hand, by  \textsc{Ekeland \& Temam} \cite[Proposition 3.1, p. 14]{EkeTem76CAVP}, there exists an affine linear minorant to $f$, i.e., there exist $\xi\in X^*$ and $\alpha\in \mathbb{R}$ such that 
\begin{align*} 
f(v)\geq \alpha+\langle \xi,v\rangle \quad \text{for all }v\in X,
\end{align*} so that $f_\varepsilon(u)>-\infty$ for every $u\in X$. This implies $\DOM(f_\varepsilon)=X$. Now, for $\lambda\in (0,1)$ and $ u_1,u_2\in X$, let $(v^i_n)_{n\in \mathbb{N}}\subset X$ be a minimizing sequence for $f_\varepsilon(u_i), i=1,2$. We set $w_n:=\lambda v_n^1+(1-\lambda)v_n^2, \, n\in \mathbb{N}$. Then, by the convexity of $f$, there holds
\begin{align*}
f_\varepsilon(\lambda u_1+(1-\lambda)u_2)&= \inf_{v\in X}\left \lbrace \frac{1}{p\varepsilon^{p-1}} \Vert \lambda u_1+(1-\lambda)u_2 -v\Vert^p+f(v)\right \rbrace\\
&\leq \frac{1}{p\varepsilon^{p-1}} \Vert \lambda u_1+(1-\lambda)u_2 -w_n\Vert^p+f(w_n)\\
&\leq \lambda\left(\frac{1}{p\varepsilon^{p-1}}\Vert u_1-v_n^1\Vert^p+f(v_n^1)\right)\\
&\quad+(1-\lambda)\left(\frac{1}{p\varepsilon^{p-1}}\Vert u_2-v_n^2\Vert^p +f(v_n^2)\right)\\
&\rightarrow\lambda f_\varepsilon(u_1)+(1-\lambda)f_\varepsilon(u_2) \quad \text{as }n\rightarrow \infty,
\end{align*} which shows the convexity of $f_\varepsilon$.  We note that by \eqref{eq:bounded}, $f_\varepsilon$ is bounded on every open bounded set of $X$. Hence, by \textsc{Ekeland \& Temam} \cite[Corollary 2.4, p. 12]{EkeTem76CAVP}, $f_\varepsilon$ is locally \textsc{Lipschitz} continuous on $X$. Finally, if $X$ is a reflexive \textsc{Banach} space, then the infimum in $f_\varepsilon(u)=\inf_{v\in X}\left\lbrace \frac{\varepsilon}{p} \left \Vert \frac{u-v}{\varepsilon}\right \Vert^p+f(v)\right\rbrace$ is attained at every point $u\in X$ by the direct method of calculus of variations.
\end{proof} In the main theorem, we will show properties of the $p$-\textsc{Moreau--Yosida} regularization under the assumption that $X$ is reflexive such that, by the renorming theorem, $X$ and $X^*$ are simultaneously strictly convex and smooth. Before we progress to the next theorem, we recall that the $p$-duality map $F_X^p$ is given by $F_X^p:=\partial \frac{1}{p}\Vert \cdot \Vert^p$ for $p>1$. Then, since the mapping $v\mapsto \frac{1}{p}\Vert v\Vert^p$ is continuous and convex on $X$, \textsc{Ekeland \& Temam} \cite[Proposition 5.1 \& 5.2, Corollary 5.1, pp. 21]{EkeTem76CAVP} ensure that $F_X^p$ is a bounded and set-valued map such that $F_X^p(u)$ is non-empty, convex, and weak*-closed for all $u\in X$. Furthermore, by \cite[Example 4.3, pp. 19]{EkeTem76CAVP} the $p$-duality map is characterized by
\begin{align}\label{eq:PDM-Caract}
F_X^p(u)=\lbrace \xi\in X^* : \langle \xi,u\rangle=\Vert u\Vert^{p}=\Vert \xi\Vert_{*}^{p^*}\rbrace.
\end{align} As for $p=2$, if the dual space is strictly convex, then by \textsc{Kien} \cite[Proposition 2.3]{BuiT02TNDM} and \textsc{Akagi \& Melchionna} \cite[Lemma 19]{AkaMel18ERNP}, the $p$-duality map is demicontinuous, single-valued, and monotone in the sense that
\begin{align*}
\langle F_X^p(u)-F_X^p(v),u-v\rangle&\geq \left( \Vert u\Vert^{p-1}-\Vert v\Vert^{p-1}\right)\left( \Vert u\Vert-\Vert v\Vert\right)\quad \text{for all }u,v\in X.
\end{align*} With the above-mentioned properties of the $p$-duality map, we are able to prove in the following theorem that the $p$-\textsc{Moreau--Yosida} regularization is, under suitable conditions, \textsc{G\^{a}teaux} differentiable and has a demicontinuous \textsc{G\^{a}teaux} derivative. This result generalizes and follows the proof of \textsc{Barbu} \cite[Theorem 2.58, p. 98]{Barb10NDMT} where the case $p=2$ has been studied. 

\begin{thm}\label{thm:MYR.Gateaux}
Let $X$ be a reflexive \textsc{Banach} space such that $X$ and its dual $X^*$ are strictly convex and smooth, and let $p>1$ and $\varepsilon>0$. Furthermore, let $f:X\rightarrow (-\infty,+\infty]$ be a proper, lower semicontinuous, and convex functional. Then, the $p$-\textsc{Moreau--Yosida} regularization is convex and locally \textsc{Lipschitz} continuous, and if $f$ is strictly convex, so is $f_\varepsilon$. Moreover, $f_\varepsilon(u)=\inf_{v\in X}\left\lbrace \frac{1}{p\varepsilon^{p-1}} \Vert u-v\Vert^p+f(v)\right\rbrace$ attains at every point $u\in X$ its unique minimizer denoted by $u_\varepsilon=J_\varepsilon(u):=\mathrm{argmin}_{v\in X}\left\lbrace \frac{1}{p\varepsilon^{p-1}} \Vert u-v\Vert^p+f(v)\right\rbrace$, and $u_\varepsilon$ satisfies the \textsc{Euler-Lagrange} equation
\begin{align}\label{eq:MYR.EL}
0\in  F^p_X\left( \frac{u_\varepsilon-u}{\varepsilon}\right)+ \partial f(u_\varepsilon).
\end{align}
Furthermore, $f_\varepsilon$ is \textsc{G\^{a}teaux}-differentiable at every point $u\in X$ with the \textsc{G\^{a}teaux}-derivative $A_\varepsilon:X\rightarrow X^*$ being demicontinuous on $X$ and satisfying $A_\varepsilon(u)=-F^p_X\left( \frac{u_\varepsilon-u}{\varepsilon}\right)$. If $X^*$ is uniformly convex, then $A_\varepsilon$ is continuous.
 Moreover, the following assertions hold:
\begin{itemize}
\item[$i)$] $f_\varepsilon(u)=\frac{\varepsilon}{p}\Vert A_\varepsilon(u)\Vert_*^{p^*}+f(u_\varepsilon)$ for every $u\in X$,
\item[$ii)$] $f(u_{\varepsilon_1})\leq f_{\varepsilon_1}(u)\leq f_{\varepsilon_2}(u)\leq f(u)$ for all $u\in X$ and all $ \varepsilon_1\geq \varepsilon_2>0$,
\item[$iii)$] $\lim_{\varepsilon\rightarrow 0}\Vert u_\varepsilon-u\Vert=0$ for all $u\in \DOM(f)$,
\item[$iv)$] $\lim_{\varepsilon \rightarrow 0}f_\varepsilon(u)=f(u)$ for every $u\in X$.
\item[$v)$] For each $u\in \DOM(\partial f)$ there holds $A_\varepsilon(u)\rightharpoonup A_0(u)\in \partial f(u)$ as $\varepsilon\rightarrow 0$, where $A_0(u):=\mathrm{argmin}\lbrace \Vert \xi\Vert_* : \xi \in \partial f(u)\rbrace$. If $X^*$ is uniformly convex, then $A_\varepsilon(u)\rightarrow A_0(u)$ as $\varepsilon\rightarrow 0$.
\end{itemize} Finally, the mapping $\varepsilon\mapsto f_\varepsilon(u)$ is differentiable on $(0,+\infty)$ with
\begin{align}\label{eq:diff.MYR}
\frac{\dd }{\dd \varepsilon}f_\varepsilon(u)=-\frac{1}{p^*\varepsilon^p}\Vert u_\varepsilon-u\Vert^p \quad \text{for all }\varepsilon>0.
\end{align}
\end{thm}
\begin{proof}[Proof] By Lemma \ref{le:MYR.basic}, the $p$-\textsc{Moreau--Yosida} regularization is convex and locally \textsc{Lipschitz} continuous on $X$. Now, let $f$ be strictly convex and let $u^0,u^1\in X$ and $t\in (0,1)$. Then, we define $u^t=tu^0+(1-t)u^1$ and assume 
\begin{align*}
f_\varepsilon(u^t)=t f_\varepsilon(u^0)+(1-t)f_\varepsilon(u^1).
\end{align*} Then, using the convexity of $\Vert \cdot\Vert^p$ and $f$, we obtain
\begin{align}\label{strict.convex}
t f_\varepsilon(u^0)+(1-t)f_\varepsilon(u^1)&=f_\varepsilon(u^t)\notag\\
&= \inf_{v\in X}\left\lbrace \frac{1}{p\varepsilon^{p-1}} \Vert u^t-v\Vert^p+f(v)\right\rbrace\notag\\
&\leq \frac{1}{p\varepsilon^{p-1}} \Vert u^t-(tu^0_\varepsilon+(1-t)u^1_\varepsilon)\Vert^p+f(tu^0_\varepsilon+(1-t)u^1_\varepsilon)\notag\\
&\leq \frac{t}{p\varepsilon^{p-1}} \Vert u^0-u^0_\varepsilon\Vert^p+\frac{(1-t)}{p\varepsilon^{p-1}} \Vert u^1-u^1_\varepsilon\Vert^p\\
&\quad+tf(u^0_\varepsilon)+(1-t)f(u^1_\varepsilon)\notag\\
&= t f_\varepsilon(u^0)+(1-t)f_\varepsilon(u^1),\notag
\end{align} where $u^i_\varepsilon:=\mathrm{argmin}_{v\in X}\left\lbrace \frac{1}{p\varepsilon^{p-1}} \Vert u^i-v\Vert^p+f(v)\right\rbrace, i=0,1$. Therefore, the inequality \eqref{strict.convex} becomes an equality that implies 
\begin{align*}
\frac{1}{p\varepsilon^{p-1}} \Vert t(u^0-u^0_\varepsilon)+(1-t)(u^1-u^1_\varepsilon)\Vert^p&=\frac{1}{p\varepsilon^{p-1}} \Vert u^t-(tu^0_\varepsilon+(1-t)u^1_\varepsilon)\Vert^p\\
&= \frac{t}{p\varepsilon^{p-1}} \Vert u^0-u^0_\varepsilon\Vert^p+\frac{(1-t)}{p\varepsilon^{p-1}} \Vert u^1-u^1_\varepsilon\Vert^p
\end{align*} and
\begin{align*}
f(tu^0_\varepsilon+(1-t)u^1_\varepsilon)&=t f(u^0_\varepsilon)+(1-t)f(u^1_\varepsilon).
\end{align*} Then, the strict convexity of the norm $\Vert \cdot \Vert$ implies $u^0-u^0_\varepsilon=u^1-u^1_\varepsilon$ and the strict convexity of $f$ implies $u^0_\varepsilon=u^1_\varepsilon$ whence $u^0=u^1$ and the strict convexity of $f_\varepsilon$.

The strict convexity of the norm also implies that the resolvent operator $J_\varepsilon(u):=\mathrm{argmin}_{v\in X}\left\lbrace \frac{1}{p\varepsilon^{p-1}} \Vert u-v\Vert^p+f(v)\right\rbrace$ is single-valued for every $u\in X$ and satisfies the inclusion \eqref{eq:MYR.EL} by \cite[Proposition 5.6, p. 26]{EkeTem76CAVP}. We define $A_\varepsilon(u):=-F^p_X\left( \frac{u_\varepsilon-u}{\varepsilon}\right)$, and note that from the characterization  \eqref{eq:PDM-Caract} of the $p$-duality map, there holds
\begin{align*}
f_\varepsilon(u)&=\frac{\varepsilon}{p}\left\Vert \frac{u_\varepsilon-u}{\varepsilon}\right\Vert^p +f(u_\varepsilon)\\
&=\frac{\varepsilon}{p}\left\Vert F_X^p\left(\frac{u_\varepsilon-u}{\varepsilon}\right)\right\Vert_*^{p^*} +f(u_\varepsilon)\\
&=\frac{\varepsilon}{p}\left\Vert A_\varepsilon(u)\right\Vert_*^{p^*} +f(u_\varepsilon).
\end{align*} If we show that the operator $A_\varepsilon$ is the \textsc{G\^{a}teaux} derivative of $f_\varepsilon$, $i)$ follows. First, \textsc{Akagi \& Melchionna} \cite[Lemma 19]{AkaMel18ERNP} have shown that the operator $A_\varepsilon:X\rightarrow X^*$ is demicontinuous, i.e., for all sequences $u_n\rightarrow u$ in $X$ as $n\rightarrow \infty$, there holds $A_\varepsilon(u_n)\rightharpoonup A_\varepsilon(u)$ in $X^*$ as $n\rightarrow \infty$. Second, we show that $A_\varepsilon(u)$ belongs to the subdifferential $\partial f_\varepsilon(u)$ for every $u\in X$. Let $u,v\in X$ and $u_\varepsilon=J_\varepsilon(u), v_\varepsilon=J_\varepsilon(v)$. Then, in view of \eqref{eq:PDM-Caract} and the fact that $A_\varepsilon(u)=-F_X^p\left( \frac{u_\varepsilon-u}{\varepsilon}\right)\in \partial f(u_\varepsilon)$,  we find
\begin{align} \label{eq:PMYR.gateaux1}
f_\varepsilon(u)-f_\varepsilon(v)&=\frac{\varepsilon}{p}\left\Vert \frac{u_\varepsilon-u}{\varepsilon}\right\Vert^p +f(u_\varepsilon)-\frac{\varepsilon}{p}\left\Vert \frac{v_\varepsilon-v}{\varepsilon}\right\Vert^p -f(v_\varepsilon)\notag \\
&\leq \frac{\varepsilon}{p}\left\Vert \frac{u_\varepsilon-u}{\varepsilon}\right\Vert^p -\frac{\varepsilon}{p}\left\Vert \frac{v_\varepsilon-v}{\varepsilon}\right\Vert^p -\left\langle F^p_X\left( \frac{u_\varepsilon-u}{\varepsilon}\right),u_\varepsilon-v_\varepsilon\right \rangle\notag\\
&= \frac{\varepsilon}{p}\left\Vert \frac{u_\varepsilon-u}{\varepsilon}\right\Vert^p -\frac{\varepsilon}{p}\left\Vert \frac{v_\varepsilon-v}{\varepsilon}\right\Vert^p -\left\langle F^p_X\left( \frac{u_\varepsilon-u}{\varepsilon}\right),u_\varepsilon-u\right \rangle \notag\\
&\quad-\left\langle F^p_X\left( \frac{u_\varepsilon-u}{\varepsilon}\right),u-v\right \rangle-\left\langle F^p_X\left( \frac{u_\varepsilon-u}{\varepsilon}\right),v-v_\varepsilon\right \rangle \notag\\
&\leq  \frac{\varepsilon}{p}\left\Vert \frac{u_\varepsilon-u}{\varepsilon}\right\Vert^p -\frac{\varepsilon}{p}\left\Vert \frac{v_\varepsilon-v}{\varepsilon}\right\Vert^p - \varepsilon\left\Vert \frac{u_\varepsilon-u}{\varepsilon}\right\Vert^p\notag \\
&\quad-\left\langle F^p_X\left( \frac{u_\varepsilon-u}{\varepsilon}\right),u-v\right \rangle+\frac{\varepsilon}{p^*}\left\Vert F^p_X \left( \frac{u_\varepsilon-u}{\varepsilon}\right) \right \Vert_*^{p^*}+\frac{\varepsilon}{p}\left\Vert\frac{v-v_\varepsilon}{\varepsilon}\right \Vert^p\notag\\
&=-\left\langle F^p_X\left( \frac{u_\varepsilon-u}{\varepsilon}\right),u-v\right \rangle \notag\\
&=\left\langle A_\varepsilon(u),u-v\right \rangle \quad \text{for all }v\in X,
\end{align} whence $A_\varepsilon(u)\in \partial f_\varepsilon(u)$. Subtracting  each side of \eqref{eq:PMYR.gateaux1} by $\left\langle A_\varepsilon(v),u-v\right \rangle$, we obtain 
\begin{align}\label{eq:gat}
0\leq f_\varepsilon(u)-f_\varepsilon(v)-\left\langle A_\varepsilon(v),u-v\right \rangle\leq \left\langle A_\varepsilon(u)-A_\varepsilon(v),u-v\right \rangle
\end{align} for all $\varepsilon>0$ and $u,v\in X$. Choosing $u=v+tw$, where $t>0$ and $w\in X$, and dividing \eqref{eq:gat} by $t$, we obtain
\begin{align*}
\lim_{t\searrow 0} \frac{f_\varepsilon(v+tw)-f_\varepsilon(v)}{t}=\left\langle A_\varepsilon(v),w\right \rangle \quad \text{ for all }w\in X,
\end{align*} where we used the demicontinuity of $A_\varepsilon$. Hence, the functional $f_\varepsilon$ is \textsc{G\^{a}teaux} differentiable with derivative $A_\varepsilon$. Adapting the proof of \cite[Proposition 1.146, p. 57]{BarPre86COBS}, we show that $A_\varepsilon$ is continuous, provided that $X^*$ is a uniformly convex space \footnote{The normed space $X$ is called uniformly convex if for each $\varepsilon\in (0,2)$ there exists $\delta(\varepsilon)>0$, for which $\Vert x\Vert\leq 1$ $\Vert y\Vert\leq 1$ and $\Vert x-y\Vert \geq \varepsilon$ imply $\left\Vert\frac{x-y}{2} \right\Vert\leq 1-\delta(\varepsilon)$.}: let $u_n\rightarrow u$ and  $u^\varepsilon_n=J_\varepsilon(u_n)$. Then, by the demicontinuity of $A_\varepsilon$, there holds 
\begin{align*}
-A_\varepsilon(u_n)= F^p_X\left( \frac{u_n^\varepsilon-u_n}{\varepsilon}\right)&\rightharpoonup F^p_X=\left( \frac{u^\varepsilon-u}{\varepsilon}\right)=-A_\varepsilon(u) \quad \text{in }X^*,\\
u_n^\varepsilon-u_n&\rightharpoonup u^\varepsilon-u\quad \text{in }X
\end{align*} as $n\rightarrow \infty$.
Then, from \eqref{eq:MYR.EL} as well as the monotonicity of the duality mapping and $\partial f$, it follows that 
\begin{align*}
0&\leq \left( \left \Vert \frac{u_n^\varepsilon-u_n}{\varepsilon}\right \Vert^{p-1}-\left \Vert \frac{u_m^\varepsilon-u_m}{\varepsilon}\right\Vert^{p-1}\right)\left( \Vert u_n^\varepsilon-u_n\Vert-\Vert u_m^\varepsilon-u_m\Vert\right)\\
&\leq \left\langle F^p_X\left( \frac{u_n^\varepsilon-u_n}{\varepsilon}\right)-F^p_X\left( \frac{u_m^\varepsilon-u_m}{\varepsilon}\right),u_n^\varepsilon-u_n-(u_m^\varepsilon-u_m)\right \rangle\leq C\Vert u_n-u_m\Vert
\end{align*} where in the last step we used the fact that $A_\varepsilon$ is demicontinuous and therefore a bounded operator. By the convergence of $(u_n)_{n\in \mathbb{N}}$, we infer that $(u^\varepsilon_n-u_n)_{n\in \mathbb{N}}$ is convergent in the norm. Since $\Vert u\Vert^p=\Vert F_X^P(u)\Vert_*^{p^*}$, the sequence $(F^p_X\left( \frac{u_n^\varepsilon-u_n}{\varepsilon}\right))_{n\in \mathbb{N}}$ also converges in the norm. Since $X^*$ is uniformly convex, norm convergence and weak convergence imply strong convergence, see, e.g., \cite[Proposition 3.32, p. 78]{Brez11FASS}, and thus the continuity of $A_\varepsilon$.

We prove now the assertion $ii)$. The chain of inequalities $f(u_\varepsilon)\leq f_\varepsilon(u)\leq f(u)$ follows immediately from the definition of the $p$-\textsc{Moreau--Yosida} regularization. To conclude $ii)$, it remains to show that the mapping $\varepsilon\mapsto f_\varepsilon(u)$ is monotonically decreasing on $(0,\infty)$ for every fixed $u\in X$. Let $u\in X$ and $0<\varepsilon_2<\varepsilon_1$. Then, by the definition of a minimizer
\begin{align}\label{eq:PMYR.Diff}
f_{\varepsilon_2}(u)&=\frac{\varepsilon_2}{p}\left \Vert \frac{u_{\varepsilon_2}-u}{\varepsilon_2}\right \Vert^p +f(u_{\varepsilon_2})\notag \\
&\leq \frac{\varepsilon_2}{p}\left \Vert \frac{u_{\varepsilon_1}-u}{\varepsilon_2}\right \Vert^p +f(u_{\varepsilon_1})\notag \\
&=\left(\frac{1}{p\varepsilon_2^{p-1}}-\frac{1}{p\varepsilon_1^{p-1}}\right)\Vert u_{\varepsilon_1}-u\Vert^p+\frac{\varepsilon_1}{p}\left \Vert \frac{u_{\varepsilon_1}-u}{\varepsilon_1}\right \Vert^p +f(u_{\varepsilon_1})\notag\\
&=\left(\frac{1}{p\varepsilon_2^{p-1}}-\frac{1}{p\varepsilon_1^{p-1}}\right)\Vert u_{\varepsilon_1}-u\Vert^p+f_{\varepsilon_1}(u)\\
&\leq f_{\varepsilon_1}(u)\notag.
\end{align} Now, we aim to show \eqref{eq:diff.MYR}. First, switching the roles of $\varepsilon_1$ and $\varepsilon_2$ in the inequality \eqref{eq:PMYR.Diff} and dividing both sides by $\varepsilon_1-\varepsilon_2>0$, we obtain the chain of inequalities
\begin{align}\label{eq:Diff.quot}
&\frac{1}{p(\varepsilon_2\varepsilon_1)^{p-1}}\left(\frac{\varepsilon_1^{p-1}-\varepsilon_2^{p-1}}{\varepsilon_1-\varepsilon_2}\right)\Vert u_{\varepsilon_2}-u\Vert^p\notag\\
&\leq -\frac{f_{\varepsilon_{1}}(u)-f_{\varepsilon_{2}}(u)}{\varepsilon_1-\varepsilon_2}\\
&\leq \frac{1}{p(\varepsilon_2\varepsilon_1)^{p-1}}\left(\frac{\varepsilon_1^{p-1}-\varepsilon_2^{p-1}}{\varepsilon_1-\varepsilon_2}\right)\Vert u_{\varepsilon_1}-u\Vert^p \notag
\end{align} for all $0<\varepsilon_2<\varepsilon_1$. Then, \eqref{eq:Diff.quot} implies
\begin{align}\label{eq:3.30}
\Vert u_{\varepsilon_2}-u\Vert\leq \Vert u_{\varepsilon_1}-u\Vert  \quad \text{for all }0<\varepsilon_2<\varepsilon_1.
\end{align} Second, since the real-valued mapping $\varepsilon\mapsto f_\varepsilon(u)$ is monotone for every fixed $u\in X$, it is, by \textsc{Lebesgue}'s differentiation theorem for monotone functions\footnote{See, e.g., \textsc{Elstrodt} \cite[Satz 4.5, p. 299]{Elst05MUIT}.}, almost everywhere differentiable and there holds  
\begin{align*}
\frac{\dd f_\varepsilon(u) }{\dd \varepsilon^+}\leq \frac{\dd f_\varepsilon(u) }{\dd \varepsilon^-} \quad \text{for all }\varepsilon>0, u\in X,
\end{align*} where $\frac{\dd f_\varepsilon(u) }{\dd \varepsilon^+}$ and $\frac{\dd f_\varepsilon(u) }{\dd \varepsilon^-}$ denote the right and left derivative of $\tilde{\varepsilon}\mapsto f_{\tilde{\varepsilon}}(u)$ in $\tilde{\varepsilon}=\varepsilon$, respectively. Let $\varepsilon>0$ and $h>0$ be sufficiently small. Then, choosing $\varepsilon_1=\varepsilon+h$ and $\varepsilon_2=\varepsilon$ in the first inequality as well as  $\varepsilon_1=\varepsilon$ and $\varepsilon_2=\varepsilon-h$ in the second inequality of \eqref{eq:Diff.quot} yields
\begin{align}\label{eq:ineq.1}
\frac{1}{p((\varepsilon+h)\varepsilon)^{p-1}}\left(\frac{(\varepsilon+h)^{p-1}-\varepsilon^{p-1}}{h}\right)\Vert u_{\varepsilon}-u\Vert^p  \leq -\frac{f_{\varepsilon+h}(u)-f_{\varepsilon}(u)}{h}
\end{align} and 
\begin{align} \label{eq:ineq.2}
-\frac{f_{\varepsilon}(u)-f_{\varepsilon-h}(u)}{h} &\leq \frac{1}{p((\varepsilon-h)\varepsilon_1)^{p-1}}\left(\frac{\varepsilon^{p-1}-(\varepsilon-h)^{p-1}}{h}\right)\Vert u_{\varepsilon-h}-u\Vert^p  \\
&\leq \frac{1}{p((\varepsilon-h)\varepsilon_1)^{p-1}}\left(\frac{\varepsilon^{p-1}-(\varepsilon-h)^{p-1}}{h}\right)\Vert u_{\varepsilon}-u\Vert^p \notag 
\end{align} respectively, where we employed inequality \eqref{eq:3.30}.
 Finally, letting $h\rightarrow 0$ in \eqref{eq:ineq.1} and \eqref{eq:ineq.2} yields
 \begin{align*}
 \frac{\dd f_\varepsilon}{\dd \varepsilon}=-\frac{1}{p^*\varepsilon^p}\Vert u_\varepsilon-u\Vert^p \quad \text{for  all }\varepsilon>0.
\end{align*} We continue with showing assertion $iii)$. Let $u\in \DOM(f)$, then the first inequality of \eqref{eq:Diff.quot} implies  
\begin{align}\label{eq:PMYR.conv}
\Vert u_{\varepsilon_2}-u\Vert^p&\leq \left(\frac{p(\varepsilon_2\varepsilon_1)^{p-1}}{\varepsilon_1^{p-1}-\varepsilon_2^{p-1}}\right)(f_{\varepsilon_{2}}(u)-f_{\varepsilon_1}(u))\\
&\leq \left(\frac{p(\varepsilon_2\varepsilon_1)^{p-1}}{\varepsilon_1^{p-1}-\varepsilon_2^{p-1}}\right)(f(u)-f_{\varepsilon_1}(u))\notag
\end{align} for all $0<\varepsilon_2<\varepsilon_1$. Thus, we obtain $\lim_{\varepsilon_2\rightarrow 0}\Vert u_{\varepsilon_2}-u\Vert= 0$. Taking into account the latter convergence and the lower semicontinuity of $f$, assertion $ii)$ yields 
\begin{align*}
f(u)&\leq \liminf_{\varepsilon\rightarrow 0}f(u_\varepsilon)\\
&\leq\liminf_{\varepsilon\rightarrow 0}f_\varepsilon(u)\\
&\leq\limsup_{\varepsilon\rightarrow 0}f_\varepsilon(u)\leq f(u) \quad \text{for all }u\in \DOM(f).
\end{align*} If $u\in X\backslash \DOM(f)$, we assume that there exists a sequence $(\varepsilon_n)_{n\in \mathbb{N}}\subset (0,\infty)$ with $\varepsilon_n\rightarrow 0$ as $n\rightarrow \infty$ such that  $f_{\varepsilon_n}(u)\leq C$ for all $n\in \mathbb{N}$ for a constant $C>0$. However, inequality \eqref{eq:PMYR.conv} yields $\lim_{n\rightarrow \infty}\Vert u_{\varepsilon_n}-u\Vert=0$, and we obtain $f(u)\leq \liminf f_{\varepsilon_n}(u)\leq C$, which is a contradiction to $u\in X\backslash \DOM(f)$. The assertion $v)$ follows the exact same lines as the proof of \cite[Proposition 1.146 $iv)$, p. 57]{BarPre86COBS}. 
\end{proof} The theorem showed us that the \textsc{Moreau--Yosida} regularization has indeed a regularizing effect. In fact, in view of assertion $iv)$ and \eqref{eq:diff.MYR}, one can interpret the \textsc{Moreau--Yosida} regularization as a regularization process described by the following \textsc{Hamilton--Jacobi} equation supplemented with an initial condition
\begin{align}\label{eq:HJE}
\begin{cases}
\frac{\partial}{\partial t}u(t,x)+\frac{1}{p^*}\Vert d_x u(t,x) \Vert^p =0, \quad x\in X, t>0\\
u(0+,x)=f(x), \quad \quad\quad \quad \quad\quad\quad x\in X,
\end{cases}
\end{align} where a solution $u:[0,\infty)\times X\rightarrow \mathbb{R}$ is given by the so-called \textsc{Lax--Oleinik} formula
\begin{align*}
u(t,x)=f_t(x)=\inf_{y\in X}\left \lbrace \frac{t}{p}\left \Vert \frac{x-y}{t}\right \Vert^p +f(y)\right \rbrace,
\end{align*} see, e.g., \textsc{Lions} \cite{Lion81SGEH}.

Moreover, we have seen to what extent these regularization and approximating properties depend on the properties of $X^*$. This, as previously mentioned, becomes clearer when $X=H$ is a \textsc{Hilbert} space. In this case, the \textsc{Moreau--Yosida} regularization is even \textsc{Fr\'{e}chet} differentiable and has a  \textsc{Lipschitz} continuous derivative with a \textsc{Lipschitz} constant equal to the reciprocal of the regularization parameter $\varepsilon$, see, e.g., \textsc{Barbu \& Precupanu} \cite[Corollary 2.59, p. 99]{BarPre86COBS}. Thanks to these nice properties of the regularization and its derivative that are only available on a \textsc{Hilbert} space, the \textsc{Moreau--Yosida} regularization is often applied to \textsc{Hilbert} spaces, see, e.g., \textsc{Bauschke \& Combettes} \cite{BauCom17CAMH} for a detailed treatise on \textsc{Hilbert} spaces. The \textsc{Moreau--Yosida} regularization is related to the so-called \textsc{Yosida} approximation, which, for a given operator $A$ and $\varepsilon>0$, refers to the operator $A_\varepsilon=\varepsilon^{-1}(I-S_\varepsilon)$, which is approximative to $A$, where $S_\varepsilon=(I+\varepsilon A)^{-1}$. The \textsc{Yosida} approximation is successfully employed in the theory of semigroups in order to generate strongly continuous semigroups as in the eminent \textsc{Hille--Yosida} theorem \cite{Hill52GSGC, Yosi48OPSG}, the nonlinear counterpart \cite{Dorr69NHYT, CraLig71GSGB}, or in the theory of maximal monotone operators in \textsc{Br\'{e}zis} \cite{Brez73OMMS}.\\\\ In the next theorem, we want to show that the $p$-\textsc{Moreau--Yosida} regularization preserves both the superlinearity or $p-$growth of a function and the \textsc{Mosco}-convergence of a sequence of functions $(f_n)_{n\in \mathbb{N}}$ for fixed $\varepsilon>0$. Furthermore, we show that the sequence $(f^{\varepsilon_n}_n)_{n\in \mathbb{N}}$ converges to $f$ as $\varepsilon_n \searrow 0$ in the sense of \textsc{Mosco}-convergence. Finally, we give an explicit formula for the \textsc{Legendre--Fenchel} transformation of the $p$-\textsc{Moreau--Yosida} regularization of a function. 
\begin{thm}\label{le:MYR.mosco} Let $f_n:X\rightarrow (-\infty,+\infty]$ be a proper, lower semicontinuous and convex functional for each $n\in \mathbb{N}$ such that 
\begin{itemize}
\item[i)] for all $N>0$, there holds 
  \begin{align*}
    \lim_{\Vert \xi\Vert_*\rightarrow +\infty}\frac{1}{\Vert
      \xi\Vert_*}\Big(\inf_{n\leq N}f^*_n(\xi)\Big)=\infty,\quad \lim_{\Vert v\Vert\rightarrow
      +\infty}\frac{1}{\Vert v\Vert}\Big(\inf_{n\leq N}f_n(v)\Big)=\infty.
  \end{align*}
\item[ii)] the sequence $f_n$ converges to $f$ in the sense of \textsc{Mosco} ($f_n \Mto f)$, i.e., for all $u\in X$
\begin{align*}
\begin{cases}
a)\quad  f(u)\leq \liminf_{n\to \infty} f_n(u_n) \quad \text{for all }u_n\rightharpoonup u \text{ in } X,\\
b)\quad \exists\, \hat{u}_n\rightarrow u \text{ in $X$ such that } f(u) \geq \limsup_{n\to \infty} f_n(\hat{u}_n).
\end{cases}
\end{align*}
\end{itemize} Furthermore, let $\varepsilon\in(0,1]$ and $p>1$. Then, the $p$-\textsc{Moreau--Yosida} regularization $f_n^\varepsilon$ satisfies $i)$ and $ii)$ and the convex conjugate of $f_n^{\varepsilon}$ is given by
\begin{align}\label{eq:MYR.cc}
f^{\varepsilon,*}_n(\xi)=\frac{\varepsilon}{p^*}\Vert \xi\Vert^{p^*}_*+f^*_n(\xi) \quad \text{for all } \xi\in X^*,n\in \mathbb{N},
\end{align} where $p^*>1$ is the conjugate exponent of $p$. Moreover, $f_n^\varepsilon$ and $f_n^{\varepsilon,*}$ are uniformly superlinear with respect to $\varepsilon>0$. Finally, for all sequences $(\varepsilon_n)_{n\in \mathbb{N}}\subset (0,1]$ with $\varepsilon_n\rightarrow 0$ as $n\rightarrow \infty$, there holds $f_{n}^{\varepsilon_n}\Mto f$.
\end{thm}
\begin{proof}[Proof] First, for each $\varepsilon>0$ and $n\in \mathbb{N}$, the regularization $f^\varepsilon_n$ is a proper, lower semicontinuous, and convex functional by Lemma \ref{le:MYR.basic}. The formula \eqref{eq:MYR.cc} follows from the calculations
\begin{align*}
f^{\varepsilon,*}_n(\xi)&=\sup_{v\in X}\left\lbrace \langle \xi,v\rangle-f^\varepsilon_n(v)\right \rbrace\\
&=\sup_{v\in X}\left\lbrace \langle \xi,v\rangle-\inf_{w\in X}\left \lbrace \frac{\varepsilon}{p}\left \Vert \frac{v-w}{\varepsilon}\right \Vert^p+f_n(w) \right \rbrace \right \rbrace\\
&=\sup_{v\in X}\sup_{w\in X}\left\lbrace \langle \xi,v\rangle-  \frac{\varepsilon}{p}\left \Vert \frac{v-w}{\varepsilon}\right \Vert^p-f_n(w) \right \rbrace \\
&=\sup_{w\in X}\sup_{v\in X}\left\lbrace \langle \xi,v\rangle-  \frac{\varepsilon}{p}\left \Vert \frac{v-w}{\varepsilon}\right \Vert^p-f_n(w) \right \rbrace \\
&=\sup_{w\in X}\left\lbrace  \sup_{v\in X}\left \lbrace \langle \xi,v-w\rangle -  \frac{\varepsilon}{p}\left \Vert \frac{v-w}{\varepsilon}\right \Vert^p \right \rbrace+\langle \xi,w\rangle-f_n(w) \right \rbrace\\
&=\sup_{w\in X}\left\lbrace  \varepsilon \sup_{v\in X}\left \lbrace \left \langle \xi,\frac{v-w}{\varepsilon}\right \rangle -  \frac{1}{p}\left \Vert \frac{v-w}{\varepsilon}\right \Vert^p \right \rbrace+\langle \xi,w\rangle-f_n(w) \right \rbrace\\
&=\sup_{w\in X}\left\lbrace  \frac{\varepsilon}{p^*}\Vert \xi\Vert^{p^*}_*+\langle \xi,w\rangle-f_n(w) \right \rbrace=\frac{\varepsilon}{p^*}\Vert \xi\Vert^{p^*}_*+f_n^*(\xi)
\end{align*} for all $\xi \in X^*$ and $u\in D$, where we have used the fact $(\frac{1}{p\varepsilon^{p-1}}\Vert \cdot \Vert^p)^*=\frac{\varepsilon}{p^*}\Vert \cdot \Vert_*^{p^*}$. The expression \eqref{eq:MYR.cc} also shows the superlinearity of $f^{\varepsilon,*}_n$ uniformly in $\varepsilon$.
We proceed by showing the superlinearity of $f^{\varepsilon}_n$. To do so, we note that the superlinearity of $f_n$ equivalently says that for all $N\in \mathbb{N}$ and $M>0$, there exists a positive real number $K>0$ such that  
\begin{align}\label{eq:coer.psi}
f_n(v)\geq M \Vert v\Vert
\end{align} for all $n\geq N$ and all $v\in X$ with $\Vert v\Vert\geq K$. The idea is to show that for the regularization $f^\varepsilon_n$, there exists for all $\tilde{N}\in \mathbb{N}$ and $\tilde{M}>0$ a positive real number $\tilde{K}>0$ independent of the parameter $\varepsilon>0$ such that \eqref{eq:coer.psi} is satisfied. So, let $\tilde{N}\in \mathbb{N}$ and $\tilde{M}>0$, then, for $N=\tilde{N}$ and $M=2\tilde{M}$, there exists $K>0$ such that \eqref{eq:coer.psi} holds. By \textsc{Young}'s inequality and the triangle inequality, we obtain  
\begin{align*}
f^\varepsilon_n(v)&=\inf_{\tilde{v}\in X}\left \lbrace \frac{1}{p\varepsilon^{p-1}}\Vert v-\tilde{v}\Vert^p+f_n(\tilde{v})\right \rbrace\\
&=\min \left \lbrace \inf_{\overset{\tilde{v}\in X}{\Vert \tilde{v}\Vert\geq K}}\left \lbrace \frac{1}{p\varepsilon^{p-1}}\Vert v-\tilde{v}\Vert^p+f_n(\tilde{v})\right \rbrace,\inf_{\overset{\tilde{v}\in X}{\Vert \tilde{v}\Vert\leq K}}\left \lbrace \frac{1}{p\varepsilon^{p-1}}\Vert v-\tilde{v}\Vert^p+f_n(\tilde{v}) \right \rbrace \right \rbrace\\
&\geq 
\min \left \lbrace \inf_{\overset{\tilde{v}\in X}{\Vert \tilde{v}\Vert\geq K}}\left \lbrace \frac{1}{p\varepsilon^{p-1}}\Vert v-\tilde{v}\Vert^p+M\Vert \tilde{v}\Vert \right \rbrace,\inf_{\overset{\tilde{v}\in X}{\Vert \tilde{v}\Vert\leq K}} \frac{1}{p\varepsilon^{p-1}}\Vert v-\tilde{v}\Vert^p  \right \rbrace\\
&\geq 
\min \left \lbrace \inf_{\overset{\tilde{v}\in X}{\Vert \tilde{v}\Vert\geq K}}\left \lbrace M\Vert v-\tilde{v}\Vert+M\Vert\tilde{v}\Vert -\frac{M^{p^*}\varepsilon}{p^*} \right \rbrace,\inf_{\overset{\tilde{v}\in X}{\Vert \tilde{v}\Vert\leq K}}\left \lbrace M\Vert v-\tilde{v}\Vert -\frac{M^{p^*}\varepsilon}{p^*} \right \rbrace \right \rbrace\\
&\geq \min \left \lbrace \left( M\Vert v\Vert- \frac{M^{p^*}}{p^*}\right), \left( M\Vert v\Vert -KM-\frac{M^{p^*}}{p^*}\right)\right \rbrace\\
&= M\Vert v\Vert -KM-\frac{M^{p^*}}{p^*}\\
&\geq \frac{M}{2}\Vert v\Vert=\tilde{M}\Vert v\Vert 
\end{align*} for all $v\in X \text{ with }\Vert v\Vert \geq \tilde{K}:=2\left( K+\frac{\tilde{M}^{p^*-1}}{p^*2^{p^*-1}}\right) \text{ and }\varepsilon\in (0,1]$. This implies the  superlinearity of $f_n^\varepsilon$ uniformly in $\varepsilon>0$, which, in turn implies the superlinearity for a fixed $\varepsilon>0$. We continue by showing that $f_n^\varepsilon$ is continuous in the sense of \textsc{Mosco}-convergence. In fact, we show that for a fixed $\varepsilon>0$, the regularization satisfies a stronger version of \textsc{Mosco}-convergence, meaning that there not only exists a recovery sequence, but that every sequence converging against the same limit is a recovery sequence. Let  $(v_n)_{n\in \mathbb{N}}\subset X$ be a weakly convergent sequence with weak limit $v\in X$. Now, let $(n_k)_{k\in \mathbb{N}}$ be a subsequence such that 
\begin{align*}
\liminf_{n\rightarrow \infty}f_{n}^\varepsilon(v_n)=\lim_{k\rightarrow \infty}f_{{n_k}}^\varepsilon(v_{n_k}).
\end{align*} For each $k\in \mathbb{N}$, we denote by $v_\varepsilon^k$ the unique minimizer of $v\mapsto \frac{1}{p\varepsilon^{p-1}}\Vert v-v_{n_k}\Vert^p+f_{{n_k}}(v)$ and note that thanks to the estimate
\begin{align}\label{eq:veps.bdd}
\frac{1}{p\varepsilon^{p-1}}\Vert v_{n_k}-v_\varepsilon^k\Vert^p\leq f{{n_k}}^\varepsilon(v_{n_k})\leq \frac{1}{p\varepsilon^{p-1}}\Vert v_{n_k} \Vert^p,
\end{align} the corresponding sequence of minimizers $(v_\varepsilon^k)_{k\in \mathbb{N}}$ is bounded. Therefore, there exists a subsequence (labelled as before) which is weakly convergent to an element $\tilde{v}_\varepsilon\in X$. Then, by the \textsc{Mosco}-convergence $f_{n}\Mto f$, we have 
\begin{align*}
f^\varepsilon(v)&\leq\frac{1}{p\varepsilon^{p-1}}\Vert v-\tilde{v}_\varepsilon\Vert^p+f(\tilde{v}_\varepsilon)\\
&\leq \liminf_{k\rightarrow \infty}\left \lbrace \frac{1}{p\varepsilon^{p-1}}\Vert v_{n_k}-v_\varepsilon^k\Vert^p+f_{{n_k}}(v_\varepsilon^k)\right \rbrace\\
&=\lim_{k\rightarrow \infty}f_{u{n_k}}^\varepsilon(v_{n_k})=\liminf_{n\rightarrow \infty}f_{n}^\varepsilon(v_n).
\end{align*} Now, let $v\in X$ be arbitrary and $(v_n)_{n\in \mathbb{N}}\subset X$ any strongly convergent sequence $v_n\rightarrow v$ as $n\rightarrow \infty$. We extract an arbitrary subsequence $(n_k)_{k\in \mathbb{N}}$, and to each $k\in \mathbb{N}$, we denote the minimizers of $v\mapsto \frac{1}{p\varepsilon^{p-1}}\Vert v-v_{n_k}\Vert^p+f_{{n_k}}(v)$ again by $v_\varepsilon^k\in X$. By $\tilde{v}_\varepsilon\in X$, we denote the weak limit of a further subsequence of the very same sequence which we labelled as before. Once more, by $ii)$, for the minimizer $v_\varepsilon$ of $f^\varepsilon(v)$, there exists a strongly convergent recovery sequence $(\hat{v}_k)_{k\in \mathbb{N}}\subset X$ such that $\hat{v}_k\rightarrow v_\varepsilon$ and $\lim_{k\rightarrow \infty}f_{{n_k}}(\hat{v}_k)=f(v_\varepsilon)$. It follows
\begin{align*}
f^\varepsilon(v)&\leq \frac{1}{p\varepsilon^{p-1}}\Vert v-\tilde{v}_\varepsilon\Vert^p+f(\tilde{v}_\varepsilon)\\
&\leq \liminf_{k\rightarrow \infty}\left \lbrace \frac{1}{p\varepsilon^{p-1}}\Vert v_{n_k}-v_\varepsilon^k\Vert^p+\Psi_{u_{n_k}}(v_\varepsilon^k)\right \rbrace\\
&= \liminf_{k\rightarrow \infty}f_{{n_k}}^\varepsilon(v_{n_k})\\
&\leq\limsup_{k\rightarrow \infty}f_{{n_k}}^\varepsilon(v_{n_k})\\
&\leq \limsup_{k\rightarrow \infty}\left \lbrace \frac{1}{p\varepsilon^{p-1}}\Vert v_{n_k}-\hat{v}_k\Vert^p+f_{{n_k}}(\hat{v}_k)\right \rbrace\\
&= \lim_{k\rightarrow \infty}\left \lbrace \frac{1}{p\varepsilon^{p-1}}\Vert v_{n_k}-\hat{v}_k\Vert^p+f_{{n_k}}(\hat{v}_k)\right \rbrace\\
 &=\frac{1}{p\varepsilon^{p-1}}\Vert v-v_\varepsilon\Vert^p+f(v_\varepsilon)=f^\varepsilon(v).
\end{align*} Therefore, every subsequence $(n_k)_{k\in \mathbb{N}}$ contains a further subsequence $(n_{k_l})_{l\in \mathbb{N}}$ such that $\lim_{l\rightarrow \infty}f_{{n_{k_l}}}^\varepsilon(v_{n_{k_l}})=f^\varepsilon(v)$. By the subsequence  principle, the convergence of the whole sequence follows. In particular, this shows $v_\varepsilon=\tilde{v}_\varepsilon$.

Finally, we show that the \textsc{Mosco}-convergence $f^{\varepsilon_n}_n \Mto f$ for all sequences of regularization parameters $(\varepsilon_n)_{n\in \mathbb{N}}\subset (0,1]$ with $\varepsilon_n\rightarrow 0$ as $n\rightarrow \infty$. As before, let the sequence $(v_n)_{n\in \mathbb{N}}\subset X$ be given such that $v_n\rightharpoonup v\in X$ as $n\rightarrow \infty $, and let $(n_k)_{k\in \mathbb{N}}$ be a subsequence such that 
\begin{align*}
\liminf_{n\rightarrow \infty}f_{n}^{\varepsilon_n}(v_n)=\lim_{k\rightarrow \infty}f_{{n_k}}^{\varepsilon_{n_k}}(v_{n_k}).
\end{align*} By $\tilde{v}_k \in X,\, k\in \mathbb{N}$, we denote again the minimizer of $f_{{n_k}}^{\varepsilon_{n_k}}(v_{n_k})$. Due to the same estimate as \eqref{eq:veps.bdd} for $(\tilde{v}_k)_{k\in \mathbb{N}}$, the sequence of minimizers is bounded and therefore sequentially compact with respect to the weak topology. So, after extracting a subsequence (labelled as before), we obtain a weak limit $\tilde{v}\in X$ such that $\tilde{v}_k\rightharpoonup \tilde{v}$ as $n\rightarrow \infty$. Now, we consider two cases:
\begin{itemize}
\item[$i)$] $\frac{1}{p\varepsilon_{n_k}^{p-1}}\Vert v_{n_k}-\tilde{v}_k\Vert^p \leq C$ for a constant $C>0$,
\item[$ii)$] $\frac{1}{p\varepsilon_{n_k}^{p-1}}\Vert v_{n_k}-\tilde{v}_k\Vert^p\rightarrow \infty$ as $k\rightarrow \infty$ after possibly extracting a further subsequence.
\end{itemize}
Ad $i)$. We immediately find $v=\tilde{v}$ and therefore $\tilde{v_k}\rightharpoonup v$ as $k\rightarrow \infty$. By the continuity of $f$ in the sense of \textsc{Mosco}-convergence, it follows
\begin{align*}
f(v)&\leq \liminf_{k\rightarrow \infty}f_{{n_k}}(\tilde{v}_{k})\\
&\leq \liminf_{k\rightarrow \infty}\left \lbrace\frac{1}{p\varepsilon_{n_k}^{p-1}}\Vert v_{n_k}-\tilde{v}_k\Vert^p +f_{{n_k}}(\tilde{v}_{k})\right \rbrace\\
&= \liminf_{k\rightarrow \infty}f_{{n_k}}^{\varepsilon_{n_k}}(v_{n_k})\\
&=\lim_{k\rightarrow \infty}f_{{n_k}}^{\varepsilon_{n_k}}(v_{n_k})\\
&= \liminf_{n\rightarrow \infty}f_{n}^{\varepsilon_n}(v_n).\\
\end{align*}
Ad $ii)$. We obtain
\begin{align*}
f(v)&\leq \lim_{k\rightarrow \infty} \left(\frac{1}{p\varepsilon_{n_k}^{p-1}}\Vert v_{n_k}-\tilde{v}_k\Vert^p\right)\\
&\leq \lim_{k\rightarrow \infty}\left \lbrace\frac{1}{p\varepsilon_{n_k}^{p-1}}\Vert v_{n_k}-\tilde{v}_k\Vert^p +f_{{n_k}}(\tilde{v}_{k})\right \rbrace\\
&=\lim_{k\rightarrow \infty}f_{{n_k}}^{\varepsilon_{n_k}}(v_{n_k})\\
&= \liminf_{n\rightarrow \infty}f_{n}^{\varepsilon_n}(v_n).\\
\end{align*} It remains to show the existence of a recovery sequence. Let $v\in X$ be arbitrarily chosen. Then, there exists a recovery sequence $(v_n)_{n\in \mathbb{N}}\subset X$ for $f$ with $v_n\rightarrow v$ as $v\rightarrow \infty$ such that $\lim_{n\rightarrow \infty}f_{n}(v_n)=f(v)$. Proceeding as before, we take an arbitrary subsequence $(n_k)_{k\in \mathbb{N}}$ and denote by $(\tilde{v}_k)_{k\in \mathbb{N}}\subset X$ again the minimizing sequence of $f_{{n_k}}^{\varepsilon_{n_k}}(v_{n_k})$. Then, we consider again the two cases $i)$ and $ii)$.\\
Ad $i)$. Since the recovery sequence is strongly convergent, it follows that $(\tilde{v}_k)_{k\in \mathbb{N}}$ is also strongly convergent with the same limit $v\in X$. We obtain
\begin{align*}
f(v)&\leq \liminf_{k\rightarrow \infty}f_{{n_k}}(\tilde{v}_{k})\\
&\leq \liminf_{k\rightarrow \infty}\left \lbrace\frac{1}{p\varepsilon_{n_k}^{p-1}}\Vert v_{n_k}-\tilde{v}_k\Vert^p +f_{{n_k}}(\tilde{v}_{k})\right \rbrace\\
&= \liminf_{k\rightarrow \infty}f_{{n_k}}^{\varepsilon_{n_k}}(v_{n_k})\\
&\leq \limsup_{k\rightarrow \infty}f_{{n_k}}^{\varepsilon_{n_k}}(v_{n_k})\\
&\leq\limsup_{k\rightarrow \infty}f_{{n_k}}(v_{n_k})\\
&= \lim_{k\rightarrow \infty}f_{{n_k}}(v_{n_k})=f(v),
\end{align*} which by the same argument as before implies the convergence of the full sequence, i.e., $\lim_{n\rightarrow \infty}f_{{n}}^{\varepsilon_{n}}(v_{n})=f(v)$.\\
Ad $ii)$. Due to $f_{n}^{\varepsilon_{n}}(v_{n})\leq f_{{n}}(v_{n}), n\in \mathbb{N}$, and the convergence of the right-hand side, this case cannot occur, which completes the proof.
\end{proof}

As mentioned above, the $p$-\textsc{Moreau--Yosida} regularization can be viewed as a regularization process described by the \textsc{Hamilton--Jacobi} equation \eqref{eq:HJE}. However, introducing the \textsc{Moreau--Yosida} regularization as a solution to the \textsc{Cauchy} problem \eqref{eq:HJE} does not seem 'natural`. Interestingly,  the regularization arises naturally when one deals with (generalized) gradient flow equations. To demonstrate this more clearly, we consider the generalized gradient flow 
\begin{align*}
-\vert u'(t)\vert^{p-2}u'(t)\in \partial E(u(t)),  \quad t>0,
\end{align*} of a functional $E:H\rightarrow (-\infty,+\infty]$ on a \textsc{Hilbert} space $H$. Discretizing the equation by the implicit \textsc{Euler} scheme leads to
\begin{align*}
-\left \vert \frac{U_\tau^n-U_\tau^{n-1}}{\tau}\right\vert^{p-2}	\frac{U_\tau^n-U_\tau^{n-1}}{\tau}\in \partial E(U_\tau^n), \quad n=1,2,\dots,N,
\end{align*} where, starting with $U_\tau^0=u_0\in \DOM(E)$, the values $U_\tau^n, n=1,\dots,N$, can under certain conditions be obtained by the variational approximation scheme
\begin{align}\label{eq:resolvent}
U_\tau^n\in J_\tau(U^{n-1}):=\mathrm{argmin}_{v\in H}\left \lbrace \frac{\tau}{p}\left \vert \frac{v-U_\tau^{n-1}}{\tau}\right\vert^p +E(v)\right \rbrace,  \quad n=1,2,\dots,N.
\end{align} 
Here, obviously the $p$-\textsc{Moreau--Yosida} regularization occurs naturally after discretizing the equation in time. The approximative values $U_\tau^n\in H$ are then defined by the $p$-\textsc{Moreau--Yosida} regularization $E_\tau$ where the regularization parameter is given by the step size $\tau$ of the time-discretization. It is also worth mentioning that the \textsc{Moreau--Yosida} regularization does not only regularize a function itself, but the associated resolvent operator $J_\tau(u)$ regularizes in a certain sense its arguments $u\in H$: the values $U_\tau^n\in \DOM(\partial E)$, which are achieved in the minimization scheme, are not only contained in the domain of the functional $E$, but also in the domain of the subdifferential $\partial E$.  The latter is also referred to as the regularizing or smoothing effect of the gradient flow equation, which means that for a given initial datum $u_0\in \DOM(E)$ (or in some cases even $u_0\in \overline{\DOM(E)}$) the solution does not only belong to the domain of $E$ but also to the domain of its subdifferential $\partial E$ for an infinitesimal larger time step, i.e, $u(t)\in \DOM(\partial E)$ for every $t>0$. It is well-known that for $p=2$ and when $E:H \rightarrow (-\infty,\infty]$ is a proper, lower semicontinuous, and convex functional, the subdifferential operator $\partial E$ is an infinitesimal generator of a $C_0$-semigroup such that $S(t)u_0=u(t)$ is the unique solution to the \textsc{Cauchy} problem
\begin{align*}
\begin{cases}
u'(t)\in -\partial E(u(t)), \quad t>0,\\
u(0)=u_0\in \overline{\DOM(E)}
\end{cases}
\end{align*} and which fulfills $S(t)u_0=\lim_{n\rightarrow \infty}J_{t/n}^n(u_0)$, where $J_{t/n}$ denotes again the resolvent operator given by \eqref{eq:resolvent}, see, e.g., \cite{Brez73OMMS, Barb76NSDE}. This property even holds true in a complete metric space under slightly weaker assumptions on the functional $E$, see \textsc{Ambrosio} et al. \cite{AmGiSa05GFMS} for a detailed discussion.


%
\bibliographystyle{my_alpha}
\bibliography{alex_pub,bib_aras}

\end{document}